# Small-scale operations on graphic sequences

Irena Rusu[*]

**Abstract**

A sequence $D = (d_1, d_2, \ldots, d_n)$ of positive integers is *graphic* if it is the degree sequence of a simple graph, called in this case a *realization* of $D$. In this paper, we introduce the operation of 2-*reduction*, that subtracts 1 from two integers of $D$ such that the resulting sequence $D'$ is graphic if and only if $D$ is graphic. We show that 2-reductions allow us to simply generate all the realizations of $D$, to prove existing characterizations of graphic sequences, as well as to propose new characterizations that highlight connections between mathematical and algorithmic aspects of graphic sequences.

**Keywords:** graphic degree sequence, Erdös-Gallai theorem, graph realization.

## 1 Introduction

A simple graph is a graph without loops or multiple edges. A sequence $D = (d_1, d_2, \ldots, d_n)$ of positive integers is *graphic* if a simple graph $G$ exists whose degree sequence is $D$. In this case, $G$ is a *realization* of $D$. Considering a given sequence $D$ as a potentially graphic sequence rises a two-fold problem: one aspect requires to decide whether $D$ is graphic or not, and the other to build a realization $G$ of a graphic sequence $D$, or even all the realizations $G$ of $D$.

The decision problem is addressed by two types of results: first, those proposing necessary and sufficient conditions (usually a set of inequalities to be verified) for a sequence to be graphic; second, those identifying operations on $D$ (which we will refer to as *graphic-stable*) yielding a new sequence $D'$ that is graphic if and only if $D$ is graphic. The former ones include the Erdös-Gallai [9], Berge [5] and Ryser [20] characterizations, to name only a few, but other results exist [21, 27, 22]. Among the graphic-stable operations, those due to Havel-Hakimi [12, 11] and Kleitman-Wang [16] are generally cited, but the one proposed by Tripathi-Thiagy [22] is another nice example.

A realization $G$ of a graphic sequence $D$ is implicitly built by the algorithms repeating the Havel-Hakimi or Kleitman-Wang operations until the resulting sequence is empty. Unfortunately, these algorithms add a potentially large set of edges to $G$ all at once, using criteria that could be described as greedy (to some extent), and do not allow to build all the realizations of a graphic sequence $D$. Some algorithms have been proposed to achieve this goal [14, 15], with quite complex theoretical aspects and running times. Many other aspects of realizations have been studied, resulting into a rich and interesting literature, that includes the identification of sequences with unique [18, 17, 23] or connected realizations [11, 24, 3], or admitting realizations containing a given subgraph [16, 19, 6, 26, 25], as well as admitting constrained [8, 10, 2] or relaxed [1] realizations.

In this paper, we propose a new graphic-stable operation called a 2-reduction, which is a small-scale operation when compared to the previous ones as it only subtracts 1 from two integers of the given sequence $D$. This corresponds to adding only one edge to a (possible) realization of $D$. The simplicity

[*]Nantes Université, École Centrale Nantes, CNRS, LS2N, UMR 6004, F-44000 Nantes, France

of this operation makes it a precise tool for building all the realizations of a graphic sequence $D$ as well as for studying the nature (graphic or not) of a sequence $D$.

The paper is organized as follows. In Section 2 we recall the main results and definitions. In Section 3 we introduce the 2-reduction operation, then we identify graphic-stable 2-reductions and we give an algorithm for constructing all the realizations of a graphic sequence $D$ in a simple way, by successively adding one edge to an initially edgeless graph. In Section 4 we show applications of 2-reductions to the characterization of graphic sequences. Section 5 is the conclusion.

## 2 Preliminaries

In the whole paper, $D = (d_1, d_2, \ldots, d_n)$ is a nonincreasing sequence of positive integers of length $|D| = n$. Given positive integers $a$ and $b$, we use the notation $[a..b]$ for the set of integers $i$ with $a \leq i \leq b$ (assuming that $a \leq b$), and the notation $a_{[b]}$ for the sequence made of the integer $a$ repeated $b$ times. The integers of $D$ are considered as *labeled* by their indices, meaning that a change in the order does not change the index of the integer (but only its position). We therefore define $p_i$ to be the position of $d_i$ in $D$. Easily, necessary conditions for $D$ to be graphic are $d_1 < n$, and even sum $\sum_{i=1}^{n} d_i$. Unless otherwise indicated, any operation on $D$ is followed by a reordering of $D$ in nonincreasing order and by the removal of all the zeros.

An index $i$ of $D$ is called *strong* if $d_i \geq i$ and *right strong* if moreover $d_i > d_{i+1}$. The *maximum strong index* of $D$ is an unavoidable parameter of $D$, defined as $m = \max\{i \,|\, d_i \geq i\}$. The *conjugate* $d_i^*$ and the *corrected conjugate* $\bar{d}_i$ of $d_i$ are defined as follows ($|S|$ is the cardinality of the set $S$):

$$d_i^* = |\{j \,|\, 1 \leq j \leq n, d_j \geq i\}|$$

$$\bar{d}_i = |\{j \,|\, 1 \leq j < i, d_j \geq i-1\}| + |\{j \,|\, i < j \leq n, d_j \geq i\}|.$$

These definitions are easily understood using Ferrers diagrams, that we present and use in Section 4. In the following, we assume that $d_i^* = 0$ for $i > d_1$ and $\bar{d}_i = 0$ for $i > d_1 + 1$. The conjugates and corrected conjugates are essential for the characterization of graphic sequences, as shown by the next theorem. Before, we need to present some inequalities, defined for an index $k$ and due to Erdös-Gallai [9], Berge [5] and Hässelbarth [13, 21] respectively. The latter inequality, proposed in [13] and proved in [21], is identical with (B$_k$) when $k \leq m$, hence the notation (BH$_k$) .

$$\sum_{i=1}^{k} d_i \leq k(k-1) + \sum_{i=k+1}^{n} \min\{d_i, k\} \tag{EG$_k$}$$

$$\sum_{i=1}^{k} d_i \leq \sum_{i=1}^{k} \bar{d}_i \tag{B$_k$}$$

$$\sum_{i=1}^{k} d_i \leq \sum_{i=1}^{k} (d_i^* - 1) \tag{BH$_k$}$$

The theorem below groups together the three characterizations of graphic sequences involving the above inequalities:

**Theorem 1** ([9, 5, 21])**.** *Let $D = (d_1, \ldots, d_n)$ be a nonincreasing sequence of positive integers with even sum $\sum_{i=1}^{n} d_i$. The following are equivalent:*

*(1) $D$ is graphic.*

*(2) $(\mathrm{EG}_k)$ holds for all $k \in [1..n]$.*

*(3) $(\mathrm{B}_k)$ holds for all $k \in [1..n]$.*

*(4) $(\mathrm{BH}_k)$ holds for all $k \in [1..m]$.*

It has been shown in [27] that it is sufficient to test $(\mathrm{EG}_k)$ only for the strong indices (*i.e.* for $k \in [1..m]$), and even only for $m$ and the set of right strong indices. Combined with the result from [18] showing that subtracting the right side of the inequality from its left side in $(\mathrm{EG}_k)$ produces the same result as in $(\mathrm{B}_k)$ when $k \in [1..m]$, the same reductions of the index set hold for $(\mathrm{B}_k)$ . With the reduction to $k \in [1..m]$, the inequalities of $(\mathrm{B}_k)$ and $(\mathrm{BH}_k)$ become identical, since $\bar{d}_i = d_i^* - 1$ when $i \in [1..m]$.

We now turn our attention towards the graphic-stable operations. The following theorem assembles two results due to Kleitman-Wang [16] (equivalence $1 \leftrightarrow 2$) and Tripathi-Tiagy [22] (equivalence $1 \leftrightarrow 3$). The particular case of the former one with $i = 1$ is the Havel-Hakimi theorem [12, 11]. The operation of *laying off* $d_i$, denoted $< i >$, consists in subtracting 1 from the $d_i$ leftmost (thus largest) integers $d_j$ with index $j \neq i$, and removing $d_i$ from $D$.

**Theorem 2** ([16, 22])**.** *Let $D = (d_1, d_2, \ldots, d_n)$ be a nonincreasing sequence of positive integers with $d_1 < n$. Let $i \in [1..n]$, and let $N > d_1$ be an integer. The following are equivalent:*

*(1) $D$ is graphic.*

*(2) the sequence $D^{<i>}$ obtained by laying off $d_i$ in $D$ is graphic.*

*(3) the sequence $D^{d_1 \to N} = (N, d_2, \ldots, d_n, 1_{[N-d_1]})$ is graphic.*

*KW-algorithm.* A simple algorithm follows to test whether a sequence $D$ is graphic or not. It performs laying off operations as long as possible and decides that $D$ is graphic if and only if the remaining sequence is empty. In parallel, a realization of $D$ is built from an edgeless graph with $n$ vertices $1, 2, \ldots, n$. For each operation $< i >$, edges are added between the vertex $i$ and each vertex $j$ whose degree is decreased by 1 when the operation $< i >$ is performed.

As indicated in the Introduction, the KW-algorithm allows one to find one or several realizations but not all of them. The absolute priority given to the largest degrees when $< i >$ is performed, as well as the necessity to add $d_i$ edges incident with $i$ at once, are greedy aspects of the KW-algorithm explaining that it cannot build all the realizations.

# 3 New results and algorithms

The 2-reduction that we propose in this section can be seen as a simple component of the laying off operation, but its precision allows a broader applicability.

Let $D = (d_1, d_2, \ldots, d_n)$ and consider $i, j \in [1..n]$ with $i \neq j$. The *2-reduction* $< i, j >$ of $D$ consists in subtracting 1 from $d_i$ and $d_j$ (regardless to their positions $p_i$ and $p_j$ in $D$, when different from their indices). In order to identify the 2-reductions that are graphic-stable operations, define $\widetilde{d_i} = d_i$ if $d_i < i$ and $\widetilde{d_i} = d_i + 1$ if $d_i \geq i$. Then $\widetilde{d_i}$ gives the position of the $d_i$-th leftmost (thus largest) integer of $D$, when $d_i$ is not counted.

**Lemma 3.1** (The 2-Reduction Lemma)**.** *Let $D = (d_1, d_2, \dots, d_n)$ be a nonincreasing sequence of positive integers with $d_1 < n$ and let $i \in [1..n]$. Consider $j \in [1..n] \setminus \{i\}$ such that:*

*(1) $d_j = d_{\widetilde{d_i}}$, or*

*(2) $d_j = d_{\widetilde{d_i}} + 1$, or*

*(3) $j > d_i + 1$, only when $d_i \geq d_k$ for all $k \in [1..n]$.*

*Then the sequence $D$ is graphic if and only if the sequence $D^{<i,j>}$ obtained by applying the 2-reduction $< i, j >$ to $D$ is graphic.*

*Proof.* Note that condition (3) implies either $i = 1$ or $d_i = d_1$ in our generic case, as $d_1$ is the largest integer, but - in general - concerns every occurrence of the largest integer $d_i$ from $D$. We consider the three cases for $j$:

(1) We show that $D^{<i>} = (D^{<i,j>})^{<i>}$ (up to some reordering). Then Theorem 2 ($1 \leftrightarrow 2$) applied twice allows us to conclude. Assume without loss of generality that $j = \widetilde{d_i}$ (otherwise replace $D$ with one of its nonincreasing reorderings having this property). Then in $D^{<i,j>}$ (rearranged in nonincreasing order) the integer $d'_j := d_j - 1$ is either on the same position $\widetilde{d_i}$ or to the right of it. We deduce that when $< i >$ is applied to $D^{<i,j>}$, this integer will not be decreased again, since $< i >$ uses the value $d'_i := d_i - 1$ and thus decreases the integers $d_1, \dots, d_{\widetilde{d_i}-1}$ (except $d_i - 1$, which is removed). In conclusion, $(D^{<i,j>})^{<i>}$ is obtained from $D$ by decreasing the integers $d_1, \dots, d_{\widetilde{d_i}}$ and removing $d_i$, hence $(D^{<i,j>})^{<i>} = D^{<i>}$.

(2) Again, $D^{<i>} = (D^{<i,j>})^{<i>}$ holds. To see it, remark that in $D^{<i,j>}$ we have $d'_j := d_j - 1 = d_{\widetilde{d_i}}$, thus a reordering of $D^{<i,j>}$ in nonincreasing order exists such that $d'_j$ is on position $\widetilde{d_i}$. Then the conclusion follows as above.

(3) We consider only the case $i = 1$. The other cases are similar after appropriately reordering $D$ in nonincreasing order. Then $\widetilde{d_1} = d_1 + 1$.

In the forward direction, let $D$ be graphic and use the KW-algorithm, by laying off $d_1$ first, to obtain a realization $G$ of $D$ that does not contain the edge $(1, j)$. Then let $w$ be adjacent to $j$ in $G$, and note that there exists at least one neighbor $v$ of 1 which is not adjacent to $w$; otherwise $d_w > d_1$, a contradiction. Then build the graph $G'$ that replaces the edges $(1, v)$ and $(j, w)$ of $G$ with the edges $(1, j)$ and $(v, w)$. As $G'$ is another realization of $D$, it is sufficient to remove the edge $(1, j)$ from $G'$ to obtain a realization of $D^{<1,j>}$.

In the backward direction, if $D^{<1,j>}$ is graphic, then use again the KW-algorithm, by laying off $d_1 - 1$ first, to obtain a realization of $D^{<1,j>}$ that does not contain the edge $(1, j)$. Then add the edge $(1, j)$ to obtain a realization of $D$.

□

Performing the 2-reduction $< i, j >$ on a sequence $D$ signifies that the edge $(i, j)$ is added to an initially edgeless (potential) realization $G$ of $D$ with $n$ vertices $1, 2, \dots, n$. We deduce the following corollary:

**Corollary 3.2.** *Let $D = (d_1, d_2, \dots, d_n)$ be a nonincreasing sequence of positive integers with $d_1 < n$. If $D$ is graphic, then for each $j \in [2..n]$ there exists a realization $G$ of $D$ that contains the edge $(1, j)$.*

**Algorithm 1:** One realization of $D$

```
Input: D, E

X ← ∅
while X ≠ E do
    let d_i be the integer in position 1 of D
    let d_j be the rightmost integer of D
        such that (i, j) ∈ E \ X
    X ← X ∪ {(i, j)}
    D ← D^<i,j>
G' ← ({1, 2, ..., n}, X)
return G'
```

**Algorithm 2:** All the realizations of $D$

```
Input: D, X, F
let d_i be such that p_i = 1
F+ = ∅
for h = n to d_i + 1 in decreasing order do
    F ← F ∪ F+
    let j be such that p_j = h
    if (i, j) ∉ X ∪ F then
        F+ = F+ ∪ {(i, j)}
        X ← X ∪ {(i, j)}
        D ← D^<i,j>
        if D is non-empty then
            Algorithm 2 (D, X, F)
        else
            H ← ({1, 2, ..., n}, X)
            return H
```

*Proof.* Let $G = (\{1, 2, \ldots, n\}, \emptyset)$. When $2 \le j \le d_1 + 1$, apply the KW-algorithm to $D$ by laying off $d_1$ first. When $j > d_1 + 1$, apply the 2-Reduction Lemma (case (3)) to obtain $D^{<1,j>}$ and thus to add the edge $(1, j)$ to $G$. Then apply the KW-algorithm to $D^{<1,j>}$ by laying off $d_1 - 1$ first. □

**Remark 3.1.** Once a 2-reduction $< i, j >$ is performed, testing whether $D^{<i,j>}$ is graphic and building a realization for $D$ are two different problems. One may use any possible criterion for testing whether $D^{<i,j>}$ is graphic, including one that - for instance - builds a realization of $D^{<i,j>}$ containing the edge $(i, j)$. In conclusion, not every realization of $D^{<i,j>}$ allows one to directly obtain a realization for $D$.

**Remark 3.2.** Laying off $d_i$ is the same as performing the 2-reduction $< i, \widetilde{d_i} >$ on $D$ and recursively on $D^{<i,\widetilde{d_i}>}$ until $d_i$ becomes zero.

The 2-reductions allow us to build all the realizations (with labeled vertices) of a graphic $D$ using a simple algorithm that successively adds edges to an initially edgeless graph, in all possible ways. Whereas other methods have been proposed for the same task [14, 15], the fact that they need to characterize and build the valid sets of edges to be added at once to a realization yield these algorithms complex. Our algorithm only makes, at each step, a difference between admitted and forbidden edges, and chooses one of the former.

We show first that any arbitrary, but known, realization $G$ of $D$ may be computed using a subset of 2-reductions from the 2-Reduction Lemma. Then we use the same 2-reductions to generate all the realizations of $D$.

**Claim 3.3.** *Let $D$ be a graphic sequence. Algorithm 1 builds any given realization $G = (\{1, 2, \ldots, n\}, E)$ of $D$ by performing $|E|$ 2-reductions $< i, j >$ from the 2-Reduction Lemma, such that $p_i = 1$ and $j \ge d_i + 1$.*

*Proof.* We apply an inductive reasoning to show that after the $k$-th execution of the *while* loop, $D$ is graphic, $X \subseteq E$ holds, $|X| = k$, and each current value $d_i$ counts the number of edges with endpoint $i$ belonging to $E \setminus X$. The case $k = 0$ is obvious.

During the $(k+1)$-th execution of the *while* loop, the $d_i$ neighbors $h$ of $i$ in $G$ such that $(i,h) \in E\backslash X$ either have their degrees on positions $2, \ldots, d_i{+}1$ of $D$, or at least one of them has its degree on a position exceeding $d_i + 1$. In all cases, when $j$ is maximum, as required by Algorithm 1, we have $j \geq d_i + 1$. The 2-Reduction Lemma thus applies, namely case (1) of it when $j = d_i + 1$ and case (3) of it when $j > d_i + 1$, since $p_i = 1$ thus $d_i$ is an occurrence of the largest integer of $D$. Hence the resulting sequence $D$ is graphic, and the other properties are easily verified, taking into account the inductive hypothesis. □

The same operations are then used in Algorithm 2 to generate all the realizations of $D$. The input of Algorithm 2 is $D$, the initially empty edge set $X$ of the graph $H = (\{1, 2, \ldots, n\}, X)$ and a set $F$, initially empty too, of forbidden edges. There are two layers in the algorithm:

- a *main layer*, devoted to the generation of the realizations, and defined by the lines that do not refer to $F^+$ (in this case $F$ is always empty).
- and a *restriction layer* that uses $F^+$ to define the permissions for each graph to evolve so as to avoid duplicate realizations. More precisely, in the *for* loop, $n - d_i$ new graphs are created from the graph $H = (\{1, 2, \ldots, n\}, X)$ in the input of the recursive call using the current sequence $D$ and the current set $F$. For each $h \in \{n, n-1, \ldots, d_1 + 1\}$ (in this order), the graph $H_h = (\{1, 2, \ldots, n\}, X \cup \{(i, j_h)\})$ is built, where $j_h$ is chosen such that $d_{j_h}$ is in position $h$ of $D$. The definition of $F^+$ allows a fixed graph $H_t$ to receive later any edge $(i, j_h)$ with $h < t$, but not an edge $(i, j_h)$ with $h > t$. Thus, when the algorithm creates $H_h$ - in decreasing order of $h$ - it inserts $(i, j_h)$ into $F^+$, which is later added to the set $F$ of permanently forbidden edges for each of the forthcoming graphs $H_t$ (which satisfy $t < h$).

**Claim 3.4.** *Let $D = (d_1, d_2, \ldots, d_n)$ be a graphic sequence in nonincreasing order. Algorithm 2 generates all the realizations of D (on labeled vertices) exactly once, and only them.*

*Proof.* We let $X_1$ and $X_2$ denote the set $X$ built in Algorithm 1 and Algorithm 2 respectively. The proof shows that, for any given realization $G = (\{1, 2, \ldots, n\}, E)$ of $D$, there exists a unique branch in the recursion tree of Algorithm 2 that builds $G$, and that along this branch the edges of $G$ are added to $X_2$ by Algorithm 2 exactly in the same order as they are added to $X_1$ by Algorithm 1.

Indeed, the first edge added to $X_1$ by Algorithm 1 is $(1, j)$ with $j \geq d_1 + 1$ maximum such that $(1, j) \in E$. Algorithm 2 also starts with $i = 1$ and begins the construction of a new graph $H$ for each $h \in \{n, n-1, \ldots, d_1 + 1\}$ in this order, thus including $h = j$ (recall that in the initial $D$, $p_j = h$ iff $j = h$). Then we deduce:

(i) *Identical construction.* The first edge added to $X_1$ by Algorithm 1 and to $X_2$ (on the appropriate branch) by Algorithm 2 is the same, namely $(1, j)$.

(ii) *When $j$ is considered, $F$ contains no edge from $G$.* Indeed, in Algorithm 2, when $j$ is considered $F$ contains only edges $(1, h)$ with $h > j$, which do not belong to $G$ by the maximality of $j$ in Algorithm 1.

(iii) *$G$ is built only once.* If an edge $(1, j')$ with $d_1 + 1 \leq j' < j$ exists in $G$, the *for* loop considers it later, when $F^+$ (and thus the current $F$) contains $(1, j)$. Thus the graph whose construction starts by the edge $(i, j')$ cannot be $G$, since the presence of $(1, j)$ in $F$ prevents $(1, j)$ from being added later to this graph.

The conclusion in the general case follows recursively, by noticing that the three affirmations above are valid when the $k$-th edge is added by both algorithms in $X_1$ and $X_2$ respectively:

(i) *Identical construction.* By the inductive hypothesis, there exists a node in the recursive tree of Algorithm 2 such that the partial executions of Algorithm 1 and Algorithm 2 (up to this node) build the same set $X_1 = X_2$ with $|X_1| = |X_2| = k-1$. Thus when the $k$-th edge must be chosen by each algorithm, the current $D$ is the same in both algorithms, and thus they consider the same $i$. Algorithm 1 then chooses the rightmost $d_j$, whereas Algorithm 2 starts a new branch for each $h$, including the value of $h$ such that $p_j = h$, assuming that $(i,j) \notin X_2 \cup F$.

(ii) *When $j$ is considered, $F$ contains no edge from $G$.* This is true by inductive hypothesis, and since the set $F^+$ contains only edges $(i, j')$ with $j'$ larger than $j$ (thus not from $G$ otherwise the maximality of $j$ in Algorithm 1 would be contradicted). Thus the condition above, requiring that $(i,j) \notin X_2 \cup F$, is verified.

(iii) *$G$ is built only once.* By the inductive hypothesis, $G$ is not built on another branch of the recursion tree. And $G$ cannot be built on two separate branches created by the current call, since when $(i,j)$ is added to the realization in progress, its insertion into $F^+$ ensures that no further call will be able to add this edge.

It remains to show that only the realizations of $D$ are output by Algorithm 2. This follows from the remark that Algorithm 2 decreases $d_k$, for $k \in [1..n]$, if and only if an edge incident to $k$ is added to the realization in progress. Thus when $d_k$ becomes 0, the number of edges incident with $k$ is $d_k$.

□

# 4 Characterizing graphic sequences

In this section, 2-reductions help us to provide (in Section 4.1) an alternate proof for the characterization of graphic sequences using the inequalities ($\mathrm{BH}_k$) for $k \in [1..m]$, as well as to propose (Section 4.2) a new characterization of graphic sequences.

Visualizing the Ferrers diagram facilitates a lot the explanations. The Ferrers diagram of $D$ represents $D$ as a sequence of $n$ columns aligned at their bases, numbered $1, 2, \ldots, n$ from left to right and made of $d_1, d_2, \ldots, d_n$ cells (or *units*) respectively. The rows are numbered $1, 2, \ldots d_1$ from bottom to top. Each unit of the Ferrers diagram is identified by its coordinates $(r, c)$, defined as its row number $r$ and its column number $c$ respectively. The *corrected Ferrers diagram* of $D$ is obtained from the Ferrers diagram by inserting an empty cell with coordinates $(i,i)$ in each column $d_i$ with $i \leq m$. Whereas each unit counts for 1, an empty cell counts for 0. Then $d^*_i$ and $\bar{d}_i$ respectively count the number of units on the $i$-th row in the Ferrers and corrected Ferrers diagram.

We need to state a few simple properties, easily deduced using Ferrers diagrams:

**Claim 4.1.** *Let $D = (d_1, \ldots, d_n)$ be a nonincreasing sequence of positive integers with $d_1 < n$. Then the following hold:*

*(a)* $d_n \leq m \leq d_1$.

*(b)* $d^*_{m+1} \leq m \leq d^*_m$.

*(c)* $\bar{d}_i$ *equals* $d^*_i - 1$ *when* $i \leq m$, $d^*_{m+1}$ *(which is equal to* $m$*) when* $i = m+1$ *and* $d^*_{i-1}$ *when* $m+1 < i \leq d_1 + 1$.

*Proof.* (a) We have $d_n \leq d_1 < n$ thus $d_{d_n} \geq d_n$, so that $d_n$ is a strong index. It follows that $d_n \leq m$. Moreover, $d_{d_1+1} \leq d_1 < d_1 + 1$, thus $d_1 + 1$ is not a strong index.

(b) If, by contradiction, $d^*_{m+1} > m$, then $d_{m+1} \geq m+1$ implying that $m+1$ is a strong index. This contradicts the definition of $m$. The other inequality directly follows from the definition of a strong index.

(c) When $i \leq m$, the row $i$ contains one empty cell in the corrected Ferrers diagram when compared to the Ferrers diagram, so the affirmation follows. The $m$ empty cells inserted in the Ferrers diagram in order to obtain the corrected Ferrers diagram "push" exactly $m$ units towards the $(m+1)$-th row, which does not contain other units since $d^*_{m+1} \leq m$ by affirmation (b). Hence $\bar{d}_{m+1} = m$. The other rows are not modified in the corrected Ferrers diagram, but are shifted by 1 towards the top, so that $\bar{d}_i = d^*_{i-1}$ when $i > m+1$. □

### 4.1 New proof for the inequalities $(\mathrm{BH}_k)$

The known graphic-stable operations (Havel-Hakimi, Kleitman-Wang), although potentially promising for proving (parts of) Theorem 1 by induction, seem to affect the given sequence $D$ too deeply to allow for even partial inheritance of any type of inequalities. Unlike these operations, the 2-reductions perform small and local changes, and make inductive reasoning possible.

We recall below Theorem 1 ($1 \leftrightarrow 3$), for which we provide a direct proof:

**Theorem 3.** *Let $D = (d_1, d_2, \ldots, d_n)$ be a nonincreasing sequence of positive integers, with even sum $s = \sum_{i=1}^n d_i$. Then $D$ is graphic if and only if the inequalities*

$$\sum_{i=1}^k d_i \leq \sum_{i=1}^k (d^*_i - 1) \qquad (\mathrm{BH}_k)$$

*are verified for $k \in [1..m]$.*

*Proof.* For both directions of the proof, the induction is on $s$ and the basic case is $D = (1,1)$, for which it is easy to check that the inequality holds for $k = m = 1$, and also that $D$ is graphic. Below, we apply 2-reductions to transform $D$ into $D'$ with smaller value $s$. The integers of $D'$ are denoted by $d'_i$ ($1 \leq i \leq |D'|$), and we define $m' = \max\{i \mid d'_i \geq i\}$, so that the inequalities $(\mathrm{BH}_k)$ applied to $D'$ are:

$$\sum_{i=1}^k d'_i \leq \sum_{i=1}^k (d^{*\prime}_i - 1) \qquad (\mathrm{BH'}_k)$$

for $k \in [1..m']$. We also assume that $d_1 > 1$ otherwise $D$ is made only of 1 and the conclusion follows easily.

*Forward direction.* Let $D$ be a nonincreasing sequence with $s \geq 4$ and $d_1 > 1$ which is graphic. Define the integer $t$ to be the largest index strictly smaller than $n$ such that $d_t = d_{d_n}$ (Note that $\widetilde{d_n} = d_n$). Then $t \geq d_n$. Consider $D' = D^{<n,t>}$, and note that in $D'$ the integers $d'_t := d_t - 1$ and $d'_n := d_n - 1$ are placed on the same positions as in $D$. Deux configurations, called A and B, may appear respectively when $d_n < d_t$ and $d_n = d_t$. Their Ferrers diagrams are drawn in Figure 1. We have $m = m'$, except in the case where $m = t = d_t = d_{d_n}$, for which $m' = m - 1$.

Since $D$ is graphic, $D'$ is graphic too by the 2-Reduction Lemma (case (1)). The inductive hypothesis implies that the inequalities $(\mathrm{BH'}_k)$ hold for $k \in [1..m']$. We show that the inequalities $(\mathrm{BH}_k)$ hold for $k \in [1..m]$.

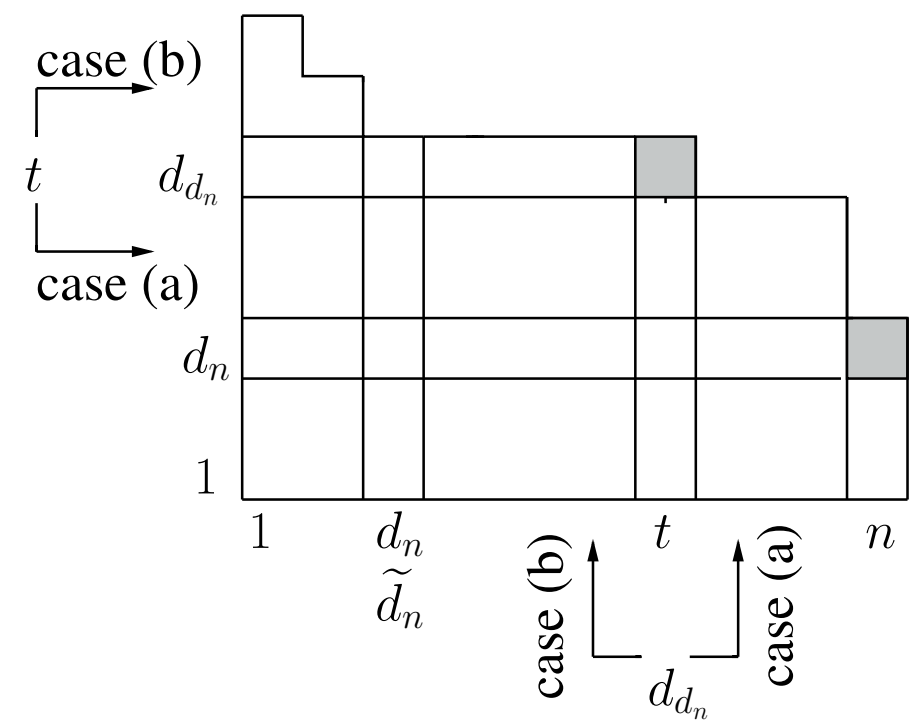

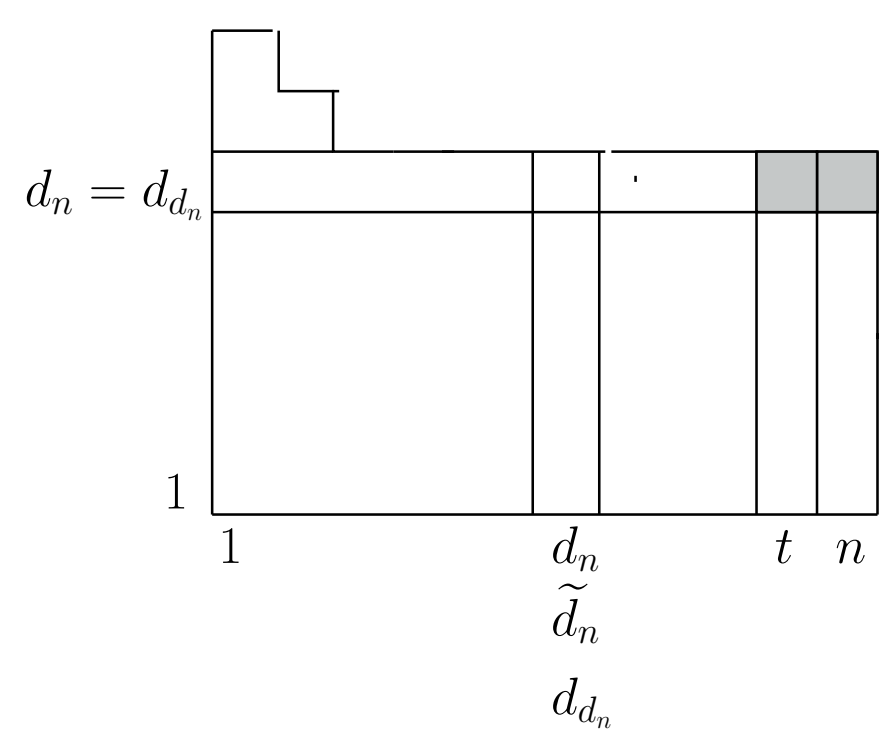

Figure 1: Ferrers diagrams of the two configurations A($d_n < d_t$) and B ($d_n = d_t$). The grey squares indicate the units removed by the 2-reduction.

**Case 1:** $m' = m$. Since $t \geq d_n$, two cases may appear:

$$\text{(a) } 1 \leq d_n \leq t \leq d_{d_n} < n \quad \text{or} \quad \text{(b) } 1 \leq d_n \leq d_{d_n} \leq t < n,$$

both with $d_n \leq m$ by Claim 4.1(a). Let $l = \sum_{i=1}^{k} d_i - \sum_{i=1}^{k} d'_i$ and $r = \sum_{i=1}^{k}(d^*_i - 1) - \sum_{i=1}^{k}(d^{*\prime}_i - 1)$, that we respectively call left and right adjustment. In order to deduce ($\mathrm{BH}_k$) from ($\mathrm{BH'}_k$), we consider all the possible cases for $k$ according to the intervals defined by (a) and by (b), and evaluate the values of $l$ and $r$. See Figure 1.

*Configuration A.* The results are presented in the *adjustment schemes* below. The intervals between two indices are close to the left and open to the right. Each of them is a possible case for $k$. The pair $(l, r)$ inside each interval $[x..y)$ indicates the values of the left and right adjustments defined above, evaluated for $k \in [x..y)$.

| | | | | | | | | | | | | | |
|---|---|---|---|---|---|---|---|---|---|---|---|---|---|
| Case (a): | 1 | [ (0,0) ) | $d_n$ | [ (0,1) ) | $t$ | [ (1,1) ) | $d_{d_n}$ | [ (1,2) ) | $n$ |
| Case (b): | 1 | [ (0,0) ) | $d_n$ | [ (0,1) ) | $d_{d_n}$ | [ (0,2) ) | $t$ | [ (1,2) ) | $n$ |

For instance, in case (a), when $k \in [d_n..t)$ we read the pair $(l, r) = (0, 1)$. The interpretation is that, for $k \in [d_n..t)$, ($\mathrm{BH}_k$) is obtained from ($\mathrm{BH'}_k$) by adding 0 to the left side of ($\mathrm{BH'}_k$) and 1 to the right side of it. The explanation is that when $k \in [d_n..t)$, $d_i = d'_i$ for all $i \leq k$, $d^*_i = d^{*\prime}_i$ for all $i \neq d_n$ but $d^*_{d_n} = d^{*\prime}_{d_n} + 1$. The other situations are explained similarly.

*Configuration B.* Only case (b) is possible, with the same adjustment scheme as above, except that the interval $[d_n..d_{d_n}]$ is empty (and possibly $[d_{d_n}..t)$ is empty too, only when $d_n = n - 1$).

In all cases $l \leq r$, so ($\mathrm{BH}_k$) is verified for all indices $k \in [1..m']$, *i.e.* for all $k \in [1..m]$.

**Case 2:** $m' = m - 1$. The inequalities for $k \leq m' = m - 1$ are proved as above, but we need to show ($\mathrm{BH}_m$). The case $m' = m - 1$ appears only when $m = t = d_t = d_{d_n}$, which implies $d_m = d_t = m = t = d_{d_n}$ and $d'_i = d_i$ for all $i < t = m$. Thus we obtain using ($\mathrm{BH'}_{m-1}$):

$$\sum_{i=1}^{m} d_i = \sum_{i=1}^{m-1} d'_i + d_m \leq \sum_{i=1}^{m-1} (d^{*\prime}_i - 1) + t \tag{1}$$

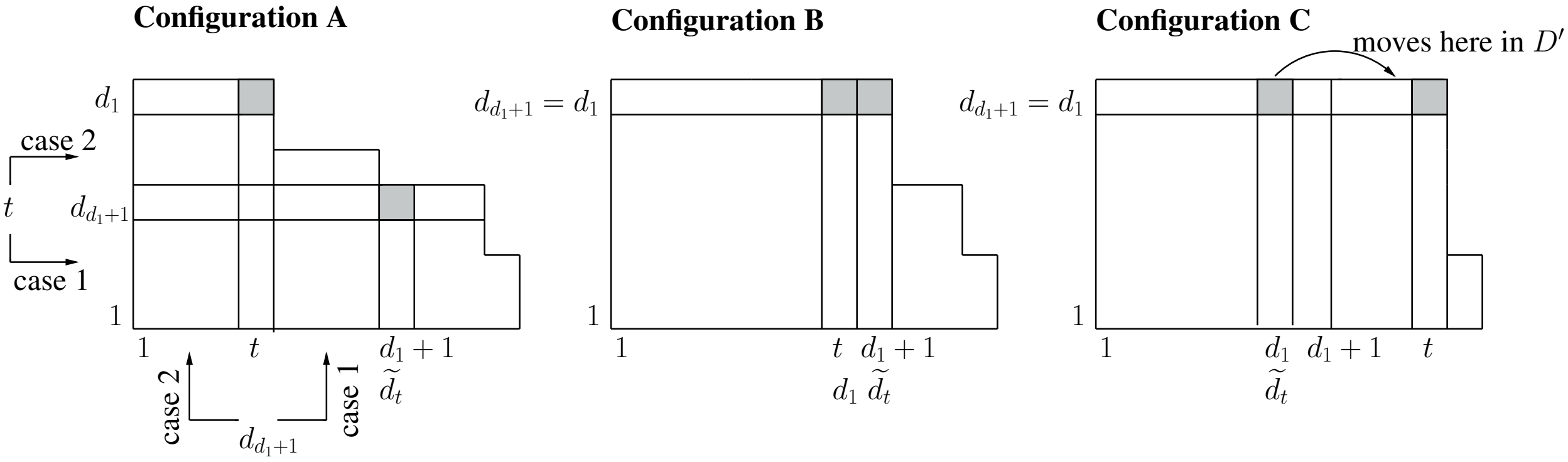

Figure 2: Ferrers diagrams of the three configurations A ($t < d_1$), B ($t = d_1$) and C ($t > d_1 + 1$). The grey squares indicate the units removed by the 2-reduction.

We also know that $d_i^{*\prime} = d_i^*$ for all $i < m$ such that $i \neq d_n$.

In configuration A, we have $d_n < d_t = m$, thus $d_{d_n}^{*\prime}$ appears in the rightmost sum of inequality (1). Since $d_{d_n}^* = d_{d_n}^{*\prime} + 1$, we deduce from (1) and from $d_m^* \geq m = t$ (use Claim 4.1(b)) that:

$$\sum_{i=1}^{m} d_i \leq \sum_{i=1}^{m-1} (d_i^* - 1) - 1 + t \leq \sum_{i=1}^{m} (d_i^* - 1),$$

which is ($\mathrm{BH}_m$). In configuration B, $d_n = d_t = t = m$ and thus (see Figure 1) we have $d_{d_n}^* = d_{d_n}^{*\prime} + 2 = t + 1$. Again, by (1) we obtain:

$$\sum_{i=1}^{m} d_i \leq \sum_{i=1}^{m-1} (d_i^{*\prime} - 1) + t \leq \sum_{i=1}^{m-1} (d_i^* - 1) + (d_{d_n}^* - 1) = \sum_{i=1}^{m} (d_i^* - 1),$$

so that ($\mathrm{BH}_m$) holds.

*Backward direction.* Let $D$ be a nonincreasing sequence with $s \geq 4$ satisfying ($\mathrm{BH}_k$) for $k \in [1..m]$. We have that $m \leq d_1$, by Claim 4.1(a), and $\widetilde{d_1} = d_1 + 1$, by definition. Let $t = d_{d_1}^*$ when $d_{d_1}^* \neq d_1 + 1$, and $t = d_1$ when $d_{d_1}^* = d_1 + 1$. Then, $\widetilde{d_t} = \widetilde{d_1} = d_1 + 1$ when $t \leq d_1$ and $\widetilde{d_t} = d_1 = \widetilde{d_1} - 1$ when $t > d_1$. Three distinct configurations, called A, B and C, exist with respect to $t$ and $d_1$, whose Ferrers diagrams are given in Figure 2. In all cases, $d_1 < n$ holds, deduced from ($\mathrm{BH}_1$). We consider $D' = D^{<t,\widetilde{d_t}>}$. In $D'$ the integer $d_t' := d_t - 1$ is still placed on position $t$ in configurations A and B, but moves in second-to-last position in configuration C. As before, $m = m'$, except in the case where $m = t = d_t = d_1$, when $m' = m - 1$. The place of the integer $d_{d_1+1} - 1$ in $D'$ is not important in configurations A and B, since it exceeds $d_1$, *i.e.* $m$.

We first show that $D'$ satisfies ($\mathrm{BH'}_k$) for $k \in [1..m']$. See Figure 2.

**Case 1:** $t \leq d_{d_1+1}$. See Configurations A (case 1) and B. We have ($\mathrm{BH}_k$) for $k \in [1..m]$ by hypothesis, $1 \leq t \leq d_{d_1+1} \leq d_1 = d_t < n$ by the properties of $t$ and $m \leq d_1$ by Claim 4.1(a). We define $l' = \sum_{i=1}^{k} d_i' - \sum_{i=1}^{k} d_i$ and $r' = \sum_{i=1}^{k} (d_i^{*\prime} - 1) - \sum_{i=1}^{k} (d_i^* - 1)$, we identify as above all the possible cases for $k$ and, for each of them, we evaluate the pair $(l', r')$. When $k < d_1$ we deduce the following adjustment scheme:

$$1 \;[\; \xrightarrow{(0,0)} \;)\; t \;[\; \xrightarrow{(-1,0)} \;)\; d_{d_1+1} \;[\; \xrightarrow{(-1,-1)} \;)\; d_1$$

Recall that the intervals are close to the left and open to the right. Since for each pair one has $l' \leq r'$, we deduce that $(\mathrm{BH'}_k)$ holds for $k < d_1$.

The case $k = d_1$ needs a proof only when $m' = m = d_1$. In this case we deduce $d_1 \geq d_m \geq m = d_1$ thus $d_m = d_1$, which implies by the definition of $t$ and from $t \leq d_{d_1+1} \leq d_1$ that $t = m = d_1$. But this is the case where $m' = m - 1$, in contradiction with the remark above that $k = d_1$ needs a proof only when $m' = m = d_1$.

**Case 2:** $d_{d_1+1} < t \leq d_1$. See configuration A (case 2). We have $(\mathrm{BH}_k)$ for $k \in [1..m]$. We also have that $t \leq m$ since $t \leq d_1 = d_t$. With the same definition for the left and right adjustments $l'$ and $r'$ as in Case 1, we deduce the following adjustment scheme:

$$1 \;[\; \xrightarrow{(0,0)} \;)\; d_{d_1+1} \;[\; \xrightarrow{(0,-1)} \;)\; t \;[\; \xrightarrow{(-1,-1)} \;)\; d_1$$

Except when $k \in [d_{d_1+1}..t)$, we remark that $l' \leq r'$ thus the inequalities $(\mathrm{BH'}_k)$ hold. When $k = d_1$ we are in the same situation as in Case 1, and we deduce a contradiction.

We therefore focus on the case $k \in [d_{d_1+1}..t)$, where the adjustment scheme shows $(l', r') = (0, -1)$. A deeper analysis is needed here. We necessarily have $n > d_1 + 1$, otherwise $(\mathrm{BH}_k)$ applied successively to $k = 1, 2, \ldots, t$ imply $d^*_1 = d^*_2 = \ldots = d^*_t = d_1 + 1$, so that $d_{d_1+1} \geq t$, a contradiction with the hypothesis of Case 2. Thus $d^*_1 = n > d_1 + 1$.

Next, the key observation is that in $D$ the lowest index $k$ for which $(\mathrm{BH}_k)$ may become an equality is $k = t$. To see it, let $h$ with $1 < h \leq t$ be the lowest index, if any, such that $\sum_{i=1}^{h} d_i = \sum_{i=1}^{h}(d^*_i - 1)$. Then, since $d^*_1 > d_1 + 1$, we deduce that $\sum_{i=2}^{h} d_i > \sum_{i=2}^{h}(d^*_i - 1)$, *i.e.* $\sum_{i=2}^{h} d_1 > \sum_{i=2}^{h}(d^*_i - 1)$. At least one of the integers $d^*_i - 1$ is then smaller than $d_1$, and - since $d^*_1, d^*_2, \ldots$ is a nonincreasing sequence - we deduce that $d^*_h - 1 < d_1$. If, by contradiction, $h < t$ holds then we have $d^*_{h+1} \leq d^*_h < d_1 + 1$ and thus:

$$\sum_{i=1}^{h+1} d_i = \sum_{i=1}^{h} d_i + d_{h+1} = \sum_{i=1}^{h}(d^*_i - 1) + d_1 > \sum_{i=1}^{h}(d^*_i - 1) + (d^*_{h+1} - 1) = \sum_{i=1}^{h+1}(d^*_i - 1).$$

This contradicts $(\mathrm{BH}_{h+1})$, which holds by hypothesis since $h + 1 \leq t \leq m$. Hence for all $k \in [d_{d_1+1}..t)$ we have the strict inequality $\sum_{i=1}^{k} d_i < \sum_{i=1}^{k}(d^*_i - 1)$, allowing us to deduce $(\mathrm{BH'}_k)$ given that in the adjustment scheme we read $(l, r) = (0, -1)$.

**Case 3:** $t > d_1 + 1$. See configuration C. We have $\widetilde{d}_t = d_1 \geq m \geq m'$. For each $k < d_1$ we deduce:

$$\sum_{i=1}^{k} d_i = kd_1 \leq k(t-2) = k(t-1) - k \leq \sum_{i=1}^{k}(d^*_i - 1) - k, \tag{2}$$

since $d_1 \leq t - 2$ given that $d_1, d_1 + 1$ and $t$ are distinct indices. The same holds for $D'$ since in the adjustment scheme we have $(l', r') = (0, 0)$. If $m' < d_1$, we are done. Otherwise, we have $m' = d_1$ and we still have to show $(\mathrm{BH'}_k)$ for $k = m' = d_1$. We have, using inequality (2) for $m' - 1$:

$$\sum_{i=1}^{m'} d'_i = \sum_{i=1}^{m'} d_i = \sum_{i=1}^{m'-1} d_i + d_1 \leq \sum_{i=1}^{m'-1} (d^*_i - 1) - (m'-1) + d_1 =$$

$$= \sum_{i=1}^{m'-1} (d^{*\prime}_i - 1) + 1 \leq \sum_{i=1}^{m'-1} (d^{*\prime}_i - 1) + (t-3) = \sum_{i=1}^{m'} (d^{*\prime}_i - 1)$$

where $d^{*\prime}_{d_1} = t - 2$ (see Figure 2) and the inequality $t \geq 4$ results from the remark that at least the distinct indices 1, $d_1$ and $d_1 + 1$ precede $t$.

By the inductive hypothesis, $D'$ is graphic, and by the 2-Reduction Lemma with $i = t$ and $j = \widetilde{d}_t$. we deduce that $D$ is graphic. □

## 4.2 New characterization of graphic sequences

The main result in this section is Theorem 4. To show it, we need some definitions and two preliminary results.

Given $D = (d_1, \ldots, d_n)$ and $i \in [1..n]$, we say that $i$ is a *weak index* of $D$ if $i \geq d_{i+2}$. Easily, if $i$ is a weak index then $i+1$ is also a weak index. For each $D$ such that $d_n < n-1$, there exists at least one weak index, namely $n-2$. The case $d_n \geq n-1$ is not particularly interesting, since then either $D$ is obviously graphic (when $d_n = n-1$, $D$ is the degree sequence of the complete graph) or obviously non-graphic (otherwise). Therefore, in the remainder of this section we assume that $d_n < n-1$. We define $w = \min\{i \mid i \geq d_{i+2}\}$, and let a *right weak index* be a weak index such that $d_j > d_{j+1}$.

The *complementary sequence* $C$ of $D$ is defined as $C = (c_1, \ldots, c_n)$ with $c_i = n - 1 - d_{n+1-i}$. It is - when $D$ is graphic - the degree sequence of the complement of any realization of $D$. We remark that $c^*_i = n - d^*_{n-i}$. The maximum strong index of $C$ is denoted by $m_C$.

**Claim 4.2.** *Let $D = (d_1, \ldots, d_n)$ be a non-increasing sequence of positive integers. The following properties hold:*

*(a)* $m - 1 \leq w \leq m$.

*(b)* $w = n - m_C - 1$.

*(c)* $w = m - 1$ *if and only if* $d^*_m = m$.

*(d)* $d^*_i \leq m$ *for all* $i \in [w+1..d_1]$.

*Proof.* (a) Since $m+1$ is not a strong index, we have $m \geq d_{m+1} \geq d_{m+2}$, implying $w \leq m$. Moreover, $d_{(m-2)+2} = d_m \geq m > m-2$, so $m-2$ is not a weak index. Thus $w \geq m-1$.

(b) We show that $m_0 := n - 1 - w$ satisfies $m_0 = m_C$. First, $m_0$ is a strong index of $C$ since $c_{m_0} = c_{n-1-w} = n - 1 - d_{n+1-(n-1-w)} = n - 1 - d_{w+2} \geq n - 1 - w = m_0$. Second, $m_0 + 1$ is not a strong index of $C$:

$$c_{m_0+1} = c_{n-1-w+1} = c_{n-w} = n - 1 - d_{n+1-(n-w)} = n - 1 - d_{w+1} < n - 1 - (w-1) = m_0 + 1,$$

where we used the definition of $w$ to deduce that $w-1$ is not a weak index, hence $w - 1 < d_{w+1}$. Thus $m_0 = m_C$.

(c) If $w = m-1$, from the property $d_{w+2} \leq w$ we obtain that $d_{m+1} \leq m-1$, thus $d^*_m \leq m$. By Claim 4.1(b) we deduce $d^*_m = m$. Conversely, $d^*_m = m$ implies that $d_{m+1} \leq m-1$, and thus $m-1$ is a weak index.

(d) From $i \geq w+1$ we deduce that $d_i^* \leq d_{w+1}^*$. If $w = m$, then $d_{w+1}^* \leq m$ by Claim 4.1(b). Otherwise, $d_{w+1}^* = d_m^* = m$ by (c). □

Before proving our theorem, we need to state the following equivalencies between inequalities, for a fixed $k$. Recall that $d_i^* = 0$ when $i > d_1$ and $\bar{d}_i = 0$ when $i > d_1 + 1$.

**Claim 4.3.** *Let $D = (d_1, \ldots, d_n)$ be a nonincreasing sequence of positive integers, with even $\sum_{i=1}^n d_i$ and $d_n < n-1$. Let $k \in [w+1..n]$. The following are equivalent:*

*(a)* $(\mathrm{B}_k)$

*(b)* $\sum_{i=k}^n d_i^* \leq \sum_{i=k+1}^n d_i$

*(c)* $\sum_{j=1}^{n-k} c_j \leq \sum_{j=1}^{n-k} \bar{c}_j$.

*Proof.* $(a) \Leftrightarrow (b)$ : By Claim 4.2(a) $k \geq w+1$ implies $k \geq m$. In the case $k > m$, by Claim 4.1(c) we deduce that:

$$\sum_{i=1}^{k} \bar{d}_i = \sum_{i=1}^{m} (d_i^* - 1) + d_{m+1}^* + \sum_{i=m+1}^{k-1} d_i^* = \sum_{i=1}^{k-1} d_i^* + d_{m+1}^* - m = \sum_{i=1}^{k-1} d_i^*.$$

The same equality follows when $k = m$: this occurs only when $w = m-1$, which is equivalent with $d_m^* = m$ by Claim 4.2(c). We deduce that:

$$\sum_{i=1}^{k} \bar{d}_i = \sum_{i=1}^{m} (d_i^* - 1) = \sum_{i=1}^{m-1} d_i^* + d_m^* - m = \sum_{i=1}^{k-1} d_i^*.$$

Then $(\mathrm{B}_k)$ is equivalent with:

$$\sum_{i=1}^{k} d_i \leq \sum_{i=1}^{k-1} d_i^*,$$

and we obtain successively (see the Ferrers diagram and note that $\sum_{i=1}^n d_i = \sum_{i=1}^n d_i^*$):

$$\sum_{i=1}^{n} d_i - \sum_{i=1}^{k-1} d_i^* \leq \sum_{i=1}^{n} d_i - \sum_{i=1}^{k} d_i$$

$$\sum_{i=k}^{n} d_i^* \leq \sum_{i=k+1}^{n} d_i.$$

$(b) \Leftrightarrow (c)$ : We start with the inequality in (b):

$$\sum_{i=k}^{n} d_i^* \leq \sum_{i=k+1}^{n} d_i.$$

As $d_n^* = 0$, $\sum_{i=k}^n d_i^* = \sum_{i=k}^{n-1} d_i^*$. Then, by a subtraction of each side from the same value $(n-k)(n-1)$ we have:

$$\sum_{i=k+1}^{n} (n-1-d_i) \leq \sum_{i=k}^{n-1} (n-1-d_i^*).$$

For the left side sum, define $j = n + 1 - i$ and note that $n - 1 - d_i = n - 1 - d_{n+1-j} = c_j$, where $k + 1 \leq n + 1 - j \leq n$, from which we deduce $1 \leq j \leq n - k$. For the right side sum, let $j = n - i$, so that $n - 1 - d_i^* = n - d_{n-j}^* - 1 = c_j^* - 1$ with $k \leq n - j \leq n - 1$, *i.e.* $1 \leq j \leq n - k$. Then the above inequality is equivalent to the one in (c), using Claim 4.1(c):

$$\sum_{j=1}^{n-k} c_j \leq \sum_{j=1}^{n-k} (c_j^* - 1) = \sum_{j=1}^{n-k} \bar{c}_j.$$

□

Consider now the TopRight-algorithm described in Algorithm 3, where $D$ is never rearranged in nonincreasing order when its integers are decreased. In the algorithm, $d_{top}$ is the integer of $D$ with largest value and smallest index. We let $right = |D|$, such that $d_{right}$ is the rightmost integer of $D$. The TopRight-algorithm successively decrements by 1 the pair of integers $d_{top}, d_{right}$ with the aim of testing whether $D$ is graphic or not. The only integers that become 0 are the ones at the end of $D$. They are systematically removed, so that $|D|$ reduces.

**Algorithm 3:** TopRight-Algorithm

**Input:** $D$
**while** $d_{top} \leq |D| - 1$ *and* $d_{top} > w$ **do**
  $D \leftarrow D^{<top,right>}$
return $(d_{top} = w)$

The next theorem first shows a link between the inequalities ($\mathrm{B}_k$) and the 2-reductions: the inequalities ($\mathrm{B}_k$) with $k \geq w+1$ are necessary and sufficient conditions for the algorithm to successively produce sequences $D$ satisfying the condition in the *while* loop, and to end with a True answer. Equivalently, and according to the proof below, the TopRight algorithm stops with $d_{top} > w$ if and only if $(\mathrm{B}_{d_{top}})$ is not verified. The algorithm may be seen as a way to test the inequalities ($\mathrm{B}_k$) for $k \geq w + 1$ and, according to the theorem below, to test if $D$ is graphic. However, as the 2-reduction is a small-scale operation, it allows precision but not necessarily rapidity. Consequently, the algorithm is not very efficient and has mainly a theoretical interest.

Secondly, the theorem shows that - similarly to the result in [27] for right strong indices - it is sufficient to test the inequalities ($\mathrm{B}_k$) for the weak indices and $w+1$ in order to deduce that $D$ is graphic.

**Theorem 4.** *Let $D = (d_1, \ldots, d_n)$ be a nonincreasing sequence of positive integers, with even $\sum_{i=1}^n d_i$ and $d_n < n - 1$. The following are equivalent:*

*(1) $D$ is graphic.*

*(2) ($\mathrm{B}_k$) holds for all $k \in [w + 1..n]$.*

*(3) the TopRight algorithm returns True.*

*(4) ($\mathrm{B}_k$) holds for $k = w + 1$ and each right weak index $k$.*

*Proof.* We show that $(1) \Rightarrow (2) \Rightarrow (4) \Rightarrow (1)$ and that $(2) \Leftrightarrow (3)$.

$(1) \Rightarrow (2)$ is a consequence of Theorem 1($1 \leftrightarrow 3$).

$(2) \Rightarrow (4)$ is obvious.

$(4) \Rightarrow (1)$ : Consider the condition in affirmation (4), that is, ($\mathrm{B}_k$) with $k = w + 1$ or $k$ a right weak index. By Claim 4.3($a \leftrightarrow c$), we have:

$$\sum_{j=1}^{n-k} c_j \leq \sum_{j=1}^{n-k} \bar{c}_j$$

for each such $k$. Let $\ell = n-k$. When $k = w+1$, $\ell = n-k = m_C$ by Claim 4.2(b) whereas when $k$ is a right weak index we have $d_{k+1} < d_k$, which is equivalent with $n-1-c_{n+1-(k+1)} < n-1-c_{n+1-k}$, that is $c_{n-k+1} < c_{n-k}$, i.e. $c_{\ell+1} < c_\ell$. Thus $\ell$ is a right strong index of $C$. But then we deduce:

$$\sum_{j=1}^{\ell} c_j \leq \sum_{j=1}^{\ell} \bar{c}_j,$$

for all $\ell$ such that either $\ell = m_C$ or $\ell$ is a right strong index of $C$.

By the result of Zverovich and Zverovich [27] (recalled in Section 2), this is sufficient to deduce that $C$ is graphic. Then $D$ is graphic too, and affirmation (1) is proved.

$(2) \Leftrightarrow (3)$ : By Claim 4.3$(a \leftrightarrow b)$, we deduce that ($\mathrm{B}_k$) for $k \in [w+1..n]$ is equivalent successively with:

$$\sum_{i=k}^{n} d_i^* \leq \sum_{i=k+1}^{n} d_i$$

$$\sum_{i=k}^{n} d_i^* - 1 < \sum_{i=k+1}^{n} d_i. \tag{3}$$

Using the Ferrers diagram, we note that the execution of the TopRight-algorithm has the following properties:

(P1) each 2-reduction $< top, right >$ in the TopRight-algorithm consists in removing the top unit from the leftmost highest column, which then represents $d_{top}$, and from the rightmost column, which represents $d_{right}$.

(P2) the rows with decreasing numbers $d_1, d_1-1, d_1-2, \ldots$ are removed in top-down order, and their units are removed in left-to-right order.

(P3) the columns with decreasing numbers $n, n-1, n-2, \ldots$ are removed in right-to-left order, and their units are removed in top-down order.

Let $k \in [w+1..n]$ be fixed, and let us consider the precise moment where the unit $(k, d_k^*)$ has to be removed by the 2-reduction $< top, right >$ in Algorithm 3. Then $d_{top} = k$ and $top = d_k^*$, and we must have $right \geq d_{top} + 1 (= k+1)$ in the condition of the *while* loop (recall that $right = |D|$). By (P1) and (P2) the number of units removed on the top of the Ferrers diagram before this 2-reduction $< top, right >$ is $\sum_{i=k}^{n} d_i^* - 1$ (all the units on the rows with number $k$ or more have been removed, except the unit $(k, d_k^*)$). By (P1) this is also the number of units removed on the right of the Ferrers diagram. Then, taking into account (P3), the inequality (3) requires that the number of units removed on the right of the Ferrers diagram, *i.e.* $\sum_{i=k}^{n} d_i^* - 1$, be strictly smaller than the number of units available in the Ferrers diagram on columns $k+1, \ldots, n$, *i.e.* $\sum_{i=k+1}^{n} d_i$. Equivalently, when $(k, d_k^*)$ is removed, there is at least one unit left on the column with number $k+1$.

This is the necessary and sufficient condition for having $|D| \geq k+1$ before the execution of the *while* loop where $d_{top} = k$ and $top = d_k^*$, but also before each execution such that $d_{top} = k$ and $top < d_k^*$. In the latter case, the value of $d_{top}$ is the same, so the condition $|D| \geq k+1$ is the same, only there are less elements already removed from the Ferrers diagram. It is immediate that when $|D| > d_{top} = k$ we also have $top \neq right$, since $k \geq w+1$ by definition, thus $top \leq d_k^* \leq m$ by Claim 4.2(d) whereas $right = |D| > d_1 \geq m$.

Thus the conditions given in inequality (3) are necessary and sufficient for the algorithm TopRight to end with $d_{top} = w$, *i.e.* to return True. □

# 5 Conclusion

We have shown in this paper that the 2-reduction we introduced is a simple tool for designing algorithms or for proving results by induction. We argue again that the simplicity of this operation allows a broad area of applications, by providing two new examples.

The first one is the proof of the Erdös-Gallai theorem by Choudum [7], that transforms $D$ into a graphically equivalent $D'$ by subtracting 1 from the rightmost largest and from the smallest integer of $D$. It is easy to notice that this very particular operation is graphic-stable. In terms of 2-reductions, we recognize $< d_1^*, n >$ from the 2-Reduction Lemma (case 3).

The second one is the following theorem due to Behzad and Chartrand [4]. A *quasi-perfect graphic sequence* is a graphic sequence whose integers are all different, except two of them which are equal to each other.

**Theorem 5.** *Let $n > 1$. There exist exactly two quasi-perfect graphic sequences of length $n$, namely:*

*(1)* $(n-1, n-2, \ldots, \lfloor n/2 \rfloor, \lfloor n/2 \rfloor, \ldots, 2, 1)$

*(2)* $(n-2, n-3, \ldots, \lfloor (n-1)/2 \rfloor, \lfloor (n-1)/2 \rfloor, \ldots, 1, 0)$

The following is an alternate proof of this result. Once the case $n = 2$ verified and the inductive hypothesis assumed, consider a quasi-perfect graphic sequence with $d_1 = n-1$ and $d_n = 1$ (as in item (1)):

$$D = (n-1, n-2, \ldots, t+1, t, t, t-1, \ldots, 2, 1).$$

We may assume that $t \leq \lfloor n/2 \rfloor$ (otherwise the same reasoning may be done on the complementary sequence). Apply successively the 2-reductions $< i, d_i >$ with $d_i = t, t-1, \ldots, 2, 1$ in this order. Then $D$ becomes $D' = (n-2, n-3, \ldots, n-t, n-(t+1), n-(t+1), n-(t+2), \ldots, t+1, t, t-1, \ldots, 1, 0)$. When the zero value is removed, $D'$ is a quasi-perfect sequence of order $n-1$ and of type (1), so by the inductive hypothesis $n-(t+1) = \lfloor (n-1)/2 \rfloor$, which implies $t = n-1-\lfloor (n-1)/2 \rfloor = \lceil (n-1)/2 \rceil = \lfloor n/2 \rfloor$. The reasoning is similar for sequences with $d_1 = n-2$ and $d_n = 0$ (as in item (2)).

We note however that, although the 2-reductions applied above do not mimic any operation of laying off, the sequence $D'$ we built above may be also obtained by laying off $d_{n-t}$ or $d_{n-(t-1)}$ (*i.e.* one of the values $t$ from $D$).

# References

[1] Amotz Bar-Noy, Toni Böhnlein, David Peleg, Mor Perry, and Dror Rawitz. Relaxed and approximate graph realizations. In *International Workshop on Combinatorial Algorithms*, pages 3–19, 2021.

[2] Amotz Bar-Noy, Igor Kalinichev, David Peleg, and Dror Rawitz. Efficient optimized degree realization: Minimum dominating set & maximum matching. *arXiv preprint arXiv:2510.03176*, 2025.

[3] Douglas Bauer, S Louis Hakimi, Nathan Kahl, and E Schmeichel. Sufficient degree conditions for $k$-edge-connectedness of a graph. *Networks: An International Journal*, 54(2):95–98, 2009.

[4] Mehdi Behzad and Gary Chartrand. No graph is perfect. *American Mathematical Monthly*, 74:962–963, 1967.

[5] Claude Berge. *Graphs and Hypergraphs*. Elsevier Science Ltd, 1985.

[6] Guantao Chen, Michael Ferrara, Ronald J Gould, and John R Schmitt. Graphic sequences with a realization containing a complete multipartite subgraph. *Discrete Mathematics*, 308(23):5712–5721, 2008.

[7] Sheshayya A Choudum. A simple proof of the erdos-gallai theorem on graph sequences. *Bulletin of the Australian Mathematical Society*, 33(1):67–70, 1986.

[8] Éva Czabarka, Aaron Dutle, Péter L Erdős, and István Miklós. On realizations of a joint degree matrix. *Discrete Applied Mathematics*, 181:283–288, 2015.

[9] Paul Erdös and Tibor Gallai. Graphs with prescribed degrees of vertices. *Mat. Lapok*, 11:264–274, 1960.

[10] Péter L Erdős, Stephen G Hartke, Leo van Iersel, and István Miklós. Graph realizations constrained by skeleton graphs. *The Electronic Journal of Combinatorics*, 24(2), 2017.

[11] S. Louis Hakimi. On realizability of a set of integers as degrees of the vertices of a linear graph. i. *Journal of the Society for Industrial and Applied Mathematics*, 10(3):496–506, 1962.

[12] Václav Havel. A remark on the existence of finite graphs. *Casopis Pest. Mat.*, 80:477–480, 1955.

[13] Werner Hässelbarth. Die verzweigtheit von graphen. *Match*, 16:3–17, 1984.

[14] Hyunju Kim, Zoltán Toroczkai, Péter L Erdős, István Miklós, and László A Székely. Degree-based graph construction. *Journal of Physics A: Mathematical and Theoretical*, 42(39):392001, 2009.

[15] Zoltán Király. Recognizing graphic degree sequences and generating all realizations. *Eötvös Loránd University, Tech. Rep. Egres TR-2011-11*, 2012.

[16] Daniel J. Kleitman and Da-Lun Wang. Algorithms for constructing graphs and digraphs with given valences and factors. *Discrete Mathematics*, 6(1):79–88, 1973.

[17] Michael Koren. Sequences with a unique realization by simple graphs. *Journal of Combinatorial Theory, Series B*, 21(3):235–244, 1976.

[18] Shuo-Yen R. Li. Graphic sequences with unique realization. *Journal of Combinatorial Theory, Series B*, 19(1):42–68, 1975.

[19] A Ramachandra Rao. The clique number of a graph with a given degree sequence. In *Proceedings of the Symposium on Graph Theory*, volume 4, pages 251–267, 1979.

[20] Herbert J Ryser. Combinatorial properties of matrices of zeros and ones. *Canadian Journal of Mathematics*, 9:371–377, 1957.

[21] Gerard Sierksma and Han Hoogeveen. Seven criteria for integer sequences being graphic. *Journal of Graph Theory*, 15(2):223–231, 1991.

[22] Amitabha Tripathi and Himanshu Tyagi. A simple criterion on degree sequences of graphs. *Discrete Applied Mathematics*, 156(18):3513–3517, 2008.

[23] Regina Tyshkevich. Decomposition of graphical sequences and unigraphs. *Discrete Mathematics*, 220(1-3):201–238, 2000.

[24] D.L. Wang and Daniel J. Kleitman. On the existence of n-connected graphs with prescribed degrees ($n \geq 2$). *Networks*, 3(3):225–239, 1973.

[25] Jian-Hua Yin. A rao-type characterization for a sequence to have a realization containing a split graph. *Discrete Mathematics*, 311(21):2485–2489, 2011.

[26] Jian-Hua Yin and Jiong-Sheng Li. Two sufficient conditions for a graphic sequence to have a realization with prescribed clique size. *Discrete Mathematics*, 301(2-3):218–227, 2005.

[27] Igor E. Zverovich and Vadim E. Zverovich. Contributions to the theory of graphic sequences. *Discrete Mathematics*, 105(1-3):293–303, 1992.